\documentclass[twoside,12pt]{article}
\usepackage{graphicx,amscd,amsfonts,amsmath,amsthm,amssymb,latexsym,amsfonts,color,extsizes,multirow,algorithm,algorithmic,mathtools}
\usepackage{epsfig,float,epstopdf,array,bm,cite,cases,array,multirow,relsize,graphicx,float,booktabs}
\usepackage[bookmarksnumbered,colorlinks,plainpages]{hyperref}
\usepackage[table]{xcolor}
\newtheorem{theorem}{\bf Theorem}
\newtheorem{lemma}{\bf Lemma}

\newtheorem{example}{Example}
\newtheorem{remark}{Remark}

\def\h{\hspace{-0.2cm}}

\begin{document}
\title{\bf A parameterized block-splitting preconditioner for  indefinite least squares problem}
\author{Davod Khojasteh Salkuyeh$^{\ddag}$\thanks {\noindent Corresponding author.}\\
	\textit{{\small $^{\ddag}$Faculty of Mathematical Sciences, University of Guilan, Rasht, Iran}} \\
	\textit{{\small  khojasteh@guilan.ac.ir, salkuyeh@gmail.com}}}
\date{}
\maketitle
\vspace{-0.5cm}

\noindent\hrulefill\\
{\bf Abstract.} We present a stationary iteration based upon a  block splitting for a class of indefinite least squares problem. Convergence of the proposed method is investigated and optimal value of the involving parameter is used. The induced preconditioner is applied to accelerate the convergence of the GMRES method for solving the problem. We also analysed the eigenpair distribution of the preconditioned matrix.  {We assess the efficiency of the proposed preconditioner by presenting results from a numerical comparison with several existing preconditioners.}   
 \\[-3mm]


\noindent\hrulefill

\pagestyle{myheadings}\markboth{D.K. Salkuyeh}{A parameterized block-splitting preconditioner for the ILS problem}
\thispagestyle{empty}

\section{Introduction}\label{s.1}
We are concerned with problem of solving the indefinite least squares (ILS)
problem 
 {
\begin{equation}\label{Eq01}
	\min_{x\in \mathbb{R}^{n}} (b-Ax)^T J (b-Ax),
\end{equation}}
where $ A\in \mathbb{R}^{m\times n} $ with $m\geq n$, and
\begin{equation}\label{Eq02}
	J=\begin{pmatrix}
		I_p & 0 \\
		0    &-I_q
	\end{pmatrix}, \qquad  p+q=m.
\end{equation}
Problems of this type are encountered in various applications, including the solution of total least squares problems \cite{BjorckBook,GolubBook,Chandrasekaran} and the field of optimization referred to as 
Robust $H^{\infty}$ smoothing \cite{Chandrasekaran}. If either $p=0$ or $q=0$ the problem is reduced to the classical least squares problem. 

 {Direct computations reveal that 
\[
f(x):=(b-Ax)^T J (b-Ax)=x^TA^TJAx-2 x^T A^T J b+b^TJb.
\]
So, we get 
\begin{equation}\label{Nabf}
\nabla f =2 (A^TJAx -A^T J b), 
\end{equation}
and the corresponding normal equations can be written as
\begin{equation}\label{Eq03}
A^T J Ax=A^TJb.
\end{equation}
It follows from  \eqref{Eq03} that  unlike the least squares problem, the ILS  problem does not necessarily have a solution. It follows from \eqref{Nabf} that the  Hessian matrix for the function $f$ is  $H=2A^T J A$.} So the minimization problem \eqref{Eq01} has a unique solution if and only if the matrix $A^T J A$ is symmetric positive definite (SPD) (see \cite{Chandrasekaran,Hassibi}).  {We assume that this condition holds throughout the paper.}

When the problem size is small, Chandrasekaran et al. in \cite{Chandrasekaran} proposed a direct method based on the QR decomposition and Cholesky factorization for solving the problem.  {Later, Bojanczyk et al. \cite{Bojanczyk1} implemented a hyperbolic QR factorization of $A$ using Householder and hyperbolic Givens transformations.
A backward stable hyperbolic QR factorization method for solving the ILS problem was presented by Xu in \cite{XuJSU}.}
However, when the problem is large, the computational cost becomes prohibitively high. So it is recommended to use iteration methods for solving the problem. 

In \cite{LiuLiu}, Liu and Liu presented the two versions of the block SOR method for solving the problem. In \cite{JSong}, Song constructed an unsymmetric SOR (USSOR) method  for solving the problem. A variable parameter Uzawa method for solving the indefinite least squares problem was presented by  Meng et al. in \cite{MengNA}. Meng et al. in \cite{MengJIAM} proposed  an enhanced Landweber method with momentum acceleration for
solving large indefinite least squares problem. Recently, Meng et al. in \cite{MengAML} have stated alternating direction implicit iteration method for the problem.

We assume that the matrix $A$ and the vector $b$ can be written as
\begin{equation}\label{Eq04}
	A=\begin{pmatrix}
		A_1 \cr
		A_2
	\end{pmatrix},\qquad 
	b=\begin{pmatrix}
		b_1 \cr
		b_2
	\end{pmatrix},
\end{equation}
where $A_1\in \mathbb{R}^{p \times n}$ is of full column rank,
$A_2\in \mathbb{R}^{q \times n}$,   $b_1\in \mathbb{R}^{p}$ and 
$b_2\in \mathbb{R}^{q}$.  Xin and Meng in \cite{Xin}   wrote the system \eqref{Eq03} in the form
\begin{equation}\label{EqBS}
	\mathcal{B} \mathbf{x}=
	\begin{pmatrix}
		I  &  A_1  & 0   \cr
		0 & P & A_2^T \cr
		0       & A_2 &I
	\end{pmatrix}
	\begin{pmatrix}
		\delta_1 \cr
		x \cr
		\delta_2
	\end{pmatrix}
	=
	\begin{pmatrix}
		b_1 \cr
		A_1^T b_1 \cr
		b_2
	\end{pmatrix},
\end{equation} 
where $I$ is the identity matrix of appropriate size, $P=A_1^TA_1$ and 
$\delta=(\delta_1;\delta_2)=b-Ax$. Then, they  considered the following three block splittings for the above system
\begin{eqnarray}
\mathcal{B} &\h=\h&	\begin{pmatrix}
	I  &  0  & 0   \cr
	0 & P & 0 \cr
	0       & 0 &I
\end{pmatrix}
-
	\begin{pmatrix}
	0  & -A_1  & 0   \cr
	0 & 0 & -A_2^T \cr
	0   & -A_2 &0
\end{pmatrix}=\mathcal{M}_1-\mathcal{N}_1, \label{BS1}
 \\
\mathcal{B} &\h=\h&	\begin{pmatrix}
	I  &  0  & 0   \cr
	0 & P & A_2^T \cr
	0       &  0 &I
\end{pmatrix}
-
\begin{pmatrix}
	0  & -A_1  & 0   \cr
	0 & 0 & 0 \cr
	0   & -A_2 &0
\end{pmatrix}=\mathcal{M}_2-\mathcal{N}_2, \label{BS2}
 \\
\mathcal{B} &\h=\h&	\begin{pmatrix}
	I  &  A_1  & 0   \cr
	0 & P & 0 \cr
	0       &  0 &I
\end{pmatrix}
-
\begin{pmatrix}
	0  & 0  & 0   \cr
	0 & 0 & -A_2^T \cr
	0   & -A_2 &0
\end{pmatrix}=\mathcal{M}_3-\mathcal{N}_3, \label{BS3}
\end{eqnarray} 
and proved that the stationary iteration methods corresponding to the above splittings are unconditionally convergent. 
 {Recently, Li et al. in \cite{LiAMC} split the matrix $\mathcal{B}$ as
\begin{equation}\label{BS4}
	\mathcal{B} =	\begin{pmatrix}
		I  &  A_1  & 0   \cr
		0  &   P   & A_2^T \cr
		0  &   0   & I
	\end{pmatrix}
	-
	\begin{pmatrix}
		0  &   0  & 0   \cr
		0  &   0  & 0 \cr
		0  & -A_2 & 0
	\end{pmatrix}=\mathcal{M}_4-\mathcal{N}_4, 
\end{equation} 
and showed that the corresponding stationary iteration method is also unconditionally convergent. The authors of \cite{LiAMC,Xin} applied the matrices  $\mathcal{M}_i$ $(i=1,2,3,4)$ as preconditioners to accelerate the convergence of GMRES for solving the system   \eqref{Eq05}. The numerical results presented in \cite{LiAMC} demonstrate that the preconditioner $\mathcal{M}_4$ outperforms the preconditioners $\mathcal{M}_i$ $(i=1,2,3)$.    }

 {
In this paper, we consider the following reformulation used in \cite{JSong}
of the problem \eqref{Eq03} 
\begin{equation}\label{Eq05}
	\begin{pmatrix}
		A_1  & 0 & I  \cr
		A_2  & I & 0 \cr
		0       & -A_2^T &A_1^T
	\end{pmatrix}
	\begin{pmatrix}
		x \cr
		\delta_2 \cr
		\delta_1
	\end{pmatrix}
	=
	\begin{pmatrix}
		b_1 \cr
		b_2 \cr
		0
	\end{pmatrix},
\end{equation} 
where  $I$ denotes the identity matrix of the appropriate dimension and $\delta=(\delta_1;\delta_2)=b-Ax$.  
By premultiplying the first row of the above system by $A_1^T$ and letting $\hat{\delta}_1=A^T\delta_1$ and $\hat{b}_1=A_1^T b$, we obtain
the following system
\begin{equation}\label{Eq06}
	\mathcal{A} \mathbf{x}\equiv \begin{pmatrix}
		P  & 0 & I  \cr
		A_2  & I & 0 \cr
		0       & -A_2^T &I
	\end{pmatrix}
	\begin{pmatrix}
		x \cr
		\delta_2 \cr
		\hat{\delta}_1
	\end{pmatrix}
	=
	\begin{pmatrix}
		\hat{b}_1 \cr
		b_2 \cr
		0
	\end{pmatrix}\equiv
	\mathbf{b},
\end{equation} 
where $P=A_1^TA_1$.  In this case, all the block diagonal matrices are nonsingular. 
It is noted that the dimension of the coefficient matrix of this system is $\textbf{n}=2n+q$.
We propose a parameterized block splitting (PBS) iteration method for the system \eqref{Eq03},   
investigate convergence of the method and analyze  the induced preconditioner.
}

Throughout this paper, $ \rho(A) $ and $ \Vert A\Vert_{2} $  denote the spectral radius and the Euclidean norm, respectively.
For a matrix $ A \in \mathbb{R}^{r\times s}$, $ A^T$  is used for the transpose of $ A $. For two vectors $ x $ and $ y $, the notation $ (x;y) $ is used for $ (x^{T} , y^{T})^{T} $. The null space of the matrix  $A\in \mathbb{R}^{m\times n}$ is defined as $\mathcal{N}(A)=\{x\in \mathbb{R}^{n}: Ax=0\}$.   

The remainder of this paper is organized as follows.  In Section \ref{Sec2} we introduce the new iteration  method  and study its convergence properties.  In Section	\ref{Sec3}, properties of the induced preconditioner are analysed. Some numerical results are given  in Section \ref{Sec4}. Finally, a short conclusion is presented in Section \ref{Sec5}.


\section {The proposed iteration method }\label{Sec2}

We now consider the stationary iteration method induced by the parametrized block splitting (PBS)
$\mathcal{A}=\mathcal{M}_{\alpha}-\mathcal{N}_{\alpha}$, where
\begin{equation}\label{Eq07}
	\mathcal{M}_{\alpha} = \begin{pmatrix}
		P  & 0 & 0 \cr
		\alpha A_2  & I & 0 \cr
		0       & -A_2^T &I
	\end{pmatrix},\quad 
\mathcal{N}_{\alpha}=\begin{pmatrix}
	0  & 0 & -I \cr
	(\alpha-1) A_2  & 0 & 0 \cr
	0       & 0 & 0
\end{pmatrix},
\end{equation} 
and $\alpha$ is a positive number. The above splitting naturally introduces the stationary iteration method
\begin{equation}\label{Eq08}
	\mathbf{x}^ {(k +1)} = \mathcal{G}_{\alpha}\mathbf{x}^{( k)} + \mathcal{M}_{\alpha}^{-1}\mathbf{b},
\end{equation}
for solving the system \eqref{Eq06}, where $ \mathbf{x}^{(0)} $  is an initial guess and the matrix 
\[
\mathcal{G}_{\alpha}=\mathcal{M}_{\alpha}^{-1} \mathcal{N}=\mathcal{I}-\mathcal{M}_{\alpha}^{-1}\mathcal{A},
\]
 is the iteration matrix. Letting 
 $\mathbf{x}^{(k)}=(x^{(k)};\delta_2^{(k)};
 \hat{\delta}_1^{(k)})$, the PBS iteration method can be rewritten as
 	\begin{align}\label{PointPro}
 	\begin{cases}
 		x^{(k+1)}=P^{-1} (-\hat{\delta}_1^{(k)} + \hat{b}_1),\\
 		\delta_2^{(k+1)}=-\alpha A_2 x^{(k+1)}+(\alpha-1)A_2 x^{(k)}+b_2,\qquad k=0,1,\ldots,\\
 		\hat{\delta}_{1}^{(k+1)}=A_2^T \delta_2^{(k+1)},
 	\end{cases}
 \end{align}
where $\mathbf{x}^{(0)}=(x^{(0)};\delta_2^{(0)};
\hat{\delta}_1^{(0)})$ is an initial guess.

 The iteration scheme  \eqref{Eq08} converges to the solution of the system \eqref{Eq06} for any arbitrary initial guess if and only if 
 $ \rho(\mathcal{G}_{\alpha})<1$ \cite{SaadBook}.
 In the sequel we investigate the convergence of the proposed method. We first recall two lemmas.
 
\begin{lemma}\label{Lem1}\cite{AxelBook}
For two real numbers $b$ and $c$, the roots of the quadratic equation $x^2-bx+c=0$ are less than one in modulus if and only if 
\[
|c|<1,\quad |b|<1+c.
\]  
\end{lemma}

\begin{lemma}\label{Lem2}
 Let $A\in \mathbb{R}^{m\times n}$ and the matrix $A^TJA$ be SPD, where the matrix $J$ was defined in Eq. \eqref{Eq02}. Suppose that the matrix $A_1$ defined in Eq. \eqref{Eq04} is of full column rank.
If $\mu$ is an eigenvalue of the matrix $Q=(A_1^TA_1)^{-1}(A_2^TA_2)$, then $0\leq \mu <1$.
\begin{proof}
Considering  that  the matrix  $A^TJA=A_1^T A_1 - A_2^T A_2$ is SPD,
it follows from  Lemma 3.2 in \cite{LiuLiu} or Theorem 7.7.3 in \cite{Horn}. 
\end{proof}
\end{lemma}

 \begin{theorem}\label{Thm1}
Assume that the conditions of Lemma \ref{Lem2} hold. Then,
 the PBS iteration method   converges to the solution of the system \eqref{Eq06} for any arbitrary initial guess, i.e., $\rho(\mathcal{G}_{\alpha})<1$, if and only if
 	\[
 	0 < \alpha < 1+\frac{1}{\mu_{\max}}. 
 	\]
 	where $\mu_{max}$ is the largest eigenvalue of  $Q=(A_1^TA_1)^{-1}(A_2^TA_2)$. 
 	\begin{proof}
 		 {Let $(\lambda,\textbf{w})$ be an eigenpair of the iteration matrix $\mathcal{G}_{\alpha}$}, i.e., $\mathcal{G}_{\alpha} \textbf{w}=\lambda \textbf{w}$. This means that
 		$\mathcal{N}_{\alpha} \textbf{w}=\lambda \mathcal{M}_{\alpha}\textbf{w}$. So if we  assume that $\textbf{w}=(x;y;z)$,  	then
\begin{align}\label{Eq09}
	\begin{cases}
		-z=\lambda Px,\\
		(\alpha-1-\alpha \lambda) A_2x=\lambda y,\\
		z=A_2^Ty.
	\end{cases}
\end{align} 
Premultiplying the second equation of \eqref{Eq09} by $A_2^T$ and using the first equation,  we get
\[
(\alpha-1-\alpha \lambda)A_2^T A_2x=\lambda A_2^Ty=\lambda z.
\] 
Now, substituting $z=-\lambda Px$ from the first equation in \eqref{Eq09} in the above equation yields
\[
(\alpha-1-\alpha \lambda)P^{-1}A_2^T A_2x=-\lambda^2 x.
\]
So, we deduce that  
\[
(\alpha-1-\alpha \lambda)\mu=-\lambda^2,
\]	
where $\mu$ is an eigenvalue of $P^{-1}A_2^T A_2$. This equation can be written as
\begin{equation}\label{Quad01}
\lambda^2-\alpha \mu \lambda+(\alpha-1)\mu=0.
\end{equation}
Now, from Lemma \ref{Lem1} we see that the method is convergent if and only if  
\[
|(\alpha-1)\mu|<1, \quad  \textrm{and} \quad \alpha \mu<1+(\alpha-1)\mu.  
\]
From Lemma \ref{Lem2}, we have $0\leq \mu <1$. For $\mu=0$, both of the above inequalities hold true. So, we assume that $\mu\neq 0$. Clearly, the second equation holds true. From the first equation we see that 
\[
0 < \alpha < 1+\frac{1}{\mu},
\]
which gives the desired results.
 \end{proof}
 \end{theorem}
 \begin{remark}
 	From Theorem \ref{Thm1} we see the $\alpha=1$ is in the convergence interval of the PBS method. Hence, the method always is convergent for $\alpha=1$.  
 \end{remark}
\begin{theorem}\label{Thm2}
Assume that the conditions of Lemma \ref{Lem2} hold. The optimal value of  the PBS method is given by
\[
\alpha_{opt}={\rm arg}\hspace{-0.6cm}\min_{\hspace{-0.2cm}0 < \alpha < 1+\frac{1}{\mu_{\max}}} \rho(\mathcal{G}_{\alpha})= \frac{2}{1+\sqrt{1-\mu_{\max}}}.
\]
and  the corresponding optimal convergence factor is 
\[
\rho(\mathcal{G}_{\alpha_{opt}})= \frac{\mu_{\max}}{1+\sqrt{1-\mu_{\max}}}.
\]
\begin{proof}
	We consider the equation \eqref{Quad01} and without lose of generality  {we assume that $\mu\neq 0$.} In this case, this equation has two roots
	\[
	\lambda_{1}=\frac{\alpha \mu}{2} + \frac{\sqrt{\Delta}}{2} \quad  \textrm{and} \quad \lambda_{2}=\frac{\alpha \mu}{2} -\frac{\sqrt{\Delta}}{2},
	\]
	where the discriminant  is given by $\Delta=\alpha^2\mu^2-4(\alpha-1)\mu$.

	We rewrite  Eq. \eqref{Quad01} in the form 
	\[
	\frac{1}{\mu} \lambda^2 = \alpha \lambda +(1-\alpha).
	\] 	
	 Now, let 
	\[
	f(\lambda)=\frac{1}{\mu} \lambda^2,\quad 
	\textrm{and} \quad g_{\alpha}(\lambda)=\alpha \lambda +(1-\alpha).
	\]
	For every $\alpha$, the line $g_{\alpha}$  passes through the point $(1,1)$, i.e.,  $g_{\alpha} (1)=1$ and the slope of $g_{\alpha} (\lambda)$ is $\alpha$. Figure \ref{FigA} depicts  the points  of intersections  of the functions $f$ and $g_{\alpha}$ for an arbitrary value of $\alpha$. Here, $\lambda_1$ and $\lambda_2$ are the roots of the quadratic equation 
	\eqref{Quad01}. This figure shows that by increasing $\alpha,$ the value of $\max\{ |\lambda_{1}|,|\lambda_{2}|\}$, decrease,
	while the line $g_{\alpha} $ gets tangent to $f(\lambda)$. In the tangent case, we have $\lambda_{1}=\lambda_{2}$ and it indicates that $\Delta=0$.
	With our assumption that $\mu\neq 0$, we obtain the quadratic equation 
	\[
	\mu\alpha^2 -4\alpha+4=0.
	\] 
	This equation has two solutions 
	\[
	\alpha_{\pm}=\frac{2}{1\pm \sqrt{1-\mu}}. 
	\]
	It is easy to see that $\alpha_{-} \geq 1+1/\mu_{\max}$, and so $\alpha_{-}$ does not belong to the convergence interval. On the other hand, for $\alpha=\alpha_+$ we have
	\[
	\eta:=|\lambda_1|=|\lambda_2 |=\frac{\alpha_+ \mu}{2}=\sqrt{(\alpha_+-1)\mu}. 
	\]
	For every $\alpha> \alpha_+$, the quadratic equation \eqref{Quad01} has two complex roots 
	\[
\lambda_{1,2}=\frac{\alpha \mu}{2} \pm i\frac{\sqrt{-\Delta}}{2}
\]
where $i$ is the imaginary unit. It is not difficult  to see that in this case 
$$|\lambda_1|=|\lambda_2 | =\sqrt{(\alpha-1)\mu} > \eta.$$ 
So, for fixed  $\mu$, the best choice for $\alpha$ is
\[
\alpha_{+}=\frac{2}{1+ \sqrt{1-\mu}}. 
\]
Summarizing the above discussion the desired results are obtained. 	
	\begin{figure}[!t]
		\centering
			\includegraphics[height=8cm,width=10cm]{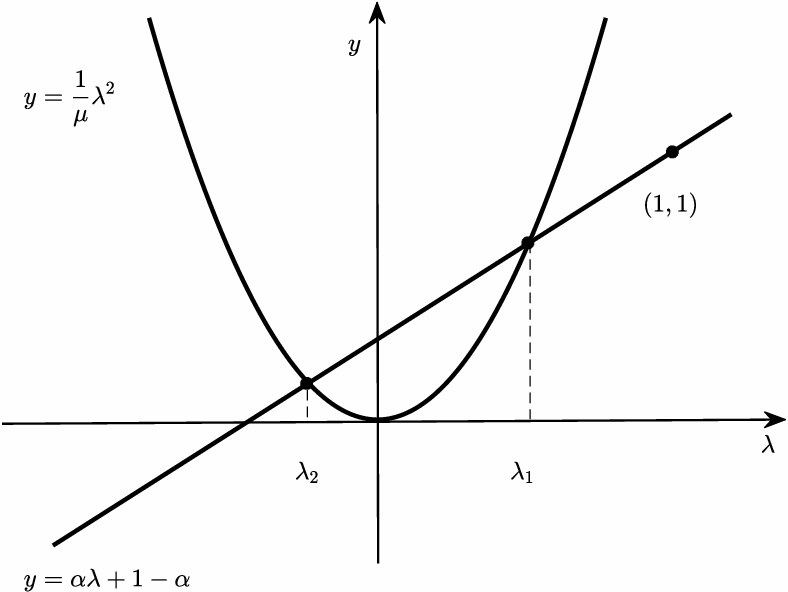}
			\caption{Depiction of the functions 
				$f(\lambda)=\frac{1}{\mu} \lambda^2$ and  $g_{\alpha}(\lambda)=\alpha \lambda +(1-\alpha)$ in Theorem \ref{Thm2}.\label{FigA}}
	\end{figure}
\end{proof}
\end{theorem}

\section{The PBS preconditioner}\label{Sec3}

In the previous section,  the PBS method is convergent if and only if  $ \alpha\in(0, 1+1/\mu_{\max})$.
It means that the  eigenvalues of  the matrix 
\begin{equation}\label{EqMNG}
\mathcal{M}_{\alpha}^{-1}\mathcal{A}=\mathcal{I}-\mathcal{G}_{\alpha},
\end{equation}
are clustered in a circle with cener $(1,0)$ and radius 1.   In this case, a Krylov subspace method  like GMRES \cite{GMRES}, would be well-suited for solving the preconditioned system \cite{BenziJCP}
\begin{equation}\label{PBSEq}
	\mathcal{M}_{\alpha}^{-1}\mathcal{A} \textbf{x} =\mathcal{M}_{\alpha}^{-1}\mathbf{b}.
\end{equation}
In each iteration of a Krylov subspace  method like GMRES it is needed to solve a system of equations
\[
 \begin{pmatrix}
	P  & 0 & 0 \cr
	\alpha A_2  & I & 0 \cr
	0       & -A_2^T &I
\end{pmatrix}
\begin{pmatrix}
	z_1 \cr
	z_2 \cr
	z_3 
\end{pmatrix}
=
\begin{pmatrix}
	r_1 \cr
	r_2 \cr
	r_3 
\end{pmatrix},
\]
can be done via the following algorithm.

\bigskip

 {\textbf{Algorithm 1. } Solution of  $\mathcal{M}_{\alpha}(z_1;z_2;z_3)=(r_1;r_2;r_3)$ \\[-0.8cm]
\begin{itemize}
	\item [1.] Solve $Pz_1=r_1$ for $z_1$;\\ [-0.8cm]
	\item [2.] Compute $z_2=r_2-\alpha A_2z_1$;\\ [-0.8cm]
	\item [3.] Compute $z_3=r_3+A_2^Tz_2$.
\end{itemize} }

The coefficient matrix of the system in Step 1 is SPD. So, it can be solved directly using the Cholesky factorization for inexactly using the conjugate gradient method.  

 {If the system of  Step 1 is solved using  the Cholesky factorization of the matrix $P$, then  the computational complexity of the algorithm is $\mathcal{O}(n^3)$, since the computational complexity of solving the system $Pz_1=r_1$ is  $\mathcal{O}(n^3)$ along with the computational complexity of matrix-vector operations in Steps 2 and 3 which is $\mathcal{O}(n^2)$. 
Comparing the PBS preconditioner with the preconditioners $\mathcal{M}_i$ ($i=1,2,3,4$)
reveals that for large systems, their computational complexities are almost the same. }

The next theorem further analyses the eigenpairs distribution of the preconditioned matrix $\mathcal{M}_{\alpha}^{-1}\mathcal{A}$. 
\begin{theorem}
Assume that the conditions of  Lemma \ref{Lem2} hold. Then, the eigenpairs of the matrix $\mathcal{M}_{\alpha}^{-1}\mathcal{A}$  are described as following
\begin{itemize}
	\item [(i)] If $\eta \neq 1$ (for every $\alpha$), then the corresponding eigenvector is given by
	\[
 \left( \frac{1}{\eta-1} P^{-1}A_2^Ty; y;A_2^Ty\right),	
	\]  
	where $y$ is nonzero vector.
	\item[(ii)] If $\alpha=\eta=1$, then every vector of the form $(x;y;0)$ with at least one of the vectors $x$ and $y$ being nonzero is the corresponding eigenvector.
	\item[(iii)] If $\alpha\neq 1$ and $\eta=1$, then every nonzero vector $(x;y;0)$ with  $x\in \mathcal{N}(A_2)$ is the corresponding eigenvector.    
\end{itemize}
\begin{proof}
 {Let $(\eta,\textbf{w})$ be an eigenpair of the matrix  $\mathcal{M}_{\alpha}^{-1}\mathcal{A}$. It follows from Eq. \eqref{EqMNG}
that each eigenvalue $\eta$ of  $\mathcal{M}_{\alpha}^{-1}\mathcal{A}$ can be expressed  as $\eta=1-\lambda$, where $\lambda$ is an eigenvalue of 
$\mathcal{G}_{\alpha}$. Substituting $\lambda=1-\eta$ in Eq. \eqref{Quad01} results in the following quadratic equation 
\[
\eta^2+(\alpha \mu-2)\eta+1-\mu=0.
\]
So each eigenvalue $\eta$ of the preconditioned   $\mathcal{M}_{\alpha}^{-1}\mathcal{A}$ is a solution of this  equation.

Letting $\textbf{w}=(x;y;z)$,  we have $\mathcal{M}_{\alpha}^{-1}\mathcal{A} \textbf{w}=\eta \textbf{w}$, which can be  written as}
\[
\begin{pmatrix}
	P  & 0 & I  \cr
	A_2  & I & 0 \cr
	0       & -A_2^T &I
\end{pmatrix}
\begin{pmatrix}
	x   \cr
	y \cr
	z
\end{pmatrix}
=
\eta \begin{pmatrix}
	P  & 0 & 0 \cr
	\alpha A_2  & I & 0 \cr
	0       & -A_2^T &I
\end{pmatrix}
\begin{pmatrix}
	x   \cr
	y \cr
	z
\end{pmatrix},
\]
which is equivalent to
	\begin{align}\label{Eq13}
		\begin{cases}
			(\eta-1)Px=z,\\
			(1-\alpha \eta) A_2x=(\eta-1) y,\\
			(\eta-1)(z-A_2^Ty)=0.
		\end{cases}
	\end{align} 
To prove $(i)$, we first show that $y\neq 0$. If $y=0$, then from the third equation in \eqref{Eq13} we get $z=0$. Substituting $z=0$ in the first equation yields $x=0$. So, in this case we deduce that $\textbf{w}=(x;y;z)=0$ which is a contradiction with the vector $\textbf{w}$ being an eigenvector. Now, from the third equation we see that $z=A_2^Ty$. Substituting $z=A_2^Ty$ in the first equation yields
\[
x=\frac{1}{\eta-1}P^{-1}z=\frac{1}{\eta-1}P^{-1}A_2^Ty.
\] 

To prove $(ii)$,  if	$\alpha=\eta=1$ then it follows from \eqref{Eq13} that $z=0$. It is noted that for every vectors $x$ and $y$, the Eq. \eqref{Eq13} holds true that proves $(ii)$.

To prove $(iii)$, we see that $z=0$ and $(1-\alpha) A_2x=0$ which implies $A_2x=0$, meaning that   $x\in \mathcal{N}(A_2)$. Hence the desired results is obtained.   
\end{proof}	
\end{theorem}  

\section{Numerical experiments}\label{Sec4}

 {In this section,} we present some numerical experiments to verify the theoretical results and show the  effectiveness of the PBS preconditioner. Comparison with other recently proposed preconditioners are   also given.  {All the numerical experiments have been conducted using some \textsc{Matlab} codes on a PC with a Pentium Dual-Core CPU, 4 gigabytes of RAM, and the Windows 7 operating system.}

\begin{example} \label{Ex1} \rm
	We consider the matrices $A_1$ and $A_2$ as following
	\[
	A_1= \begin{pmatrix}
6  &   1  &   1 \\
2  &  4   &  5  \\
1  &  1   &  5
	\end{pmatrix},
\quad 
A_2= \begin{pmatrix}	
2  &   1  &   1 \\
1  &   1  &   1 \\
1  &   2  &   2 \\
0   &  1   &  1
	\end{pmatrix},
\]
and  $b_1$ and $b_2$ are  two vectors whose entries are all ones. We have 
\[
A^TJA=A_1^TA_1-A_2^TA_2= \begin{pmatrix}	
 35 &   10  &  16 \\
10  &  11  &  19  \\
16  &  19 &   44  \\
\end{pmatrix},
\]
 {which is a SPD matrix}. So all the conditions of Theorems \ref{Thm1} and \ref{Thm2} hold true. On the other hand,
 we have $\mu_{\max}= 0.4976$.  So, the method is convergent for any $\alpha$ with
 \[
 0 < \alpha < 1+\frac{1}{\mu_{\max}}=1+\frac{1}{0.4976}= 3.009.
 \]  
 and optimal value for $\alpha$ is 
\[
\alpha_{opt}=\frac{2}{1+\sqrt{1-0.4976}}= 1.1704.
\]
Also, we have
\[
\rho(\mathcal{G}_{\alpha_{opt}})= \frac{0.4976}{1+\sqrt{1-0.4976}}=0.2912.
\]
We have solved the system  \eqref{Eq06} using the PBS method (Eq. \eqref{Eq09}). The initial guess is set to be a zeros vector and the iteration is stopped as soon as residual norm of the system is reduced by a factor of $10^{11}$. The number of iterations for the convergence and the CPU 
CPU processing time in second  for the optimum value of  $\alpha$ and some other values have been presented in Table \ref{Tbl1}. The numerical results confirm the theoretical results. 
\begin{table}[!htp]
	\centering
	\caption{Numerical results for Example \ref{Ex1}. \label{Tbl1}}\vspace{0.25cm}
	\begin{tabular}{|c|c|c|c|c|c|c|c| }
		\hline
	$\alpha$   &  ~~0.7~~    &  ~~ 0.8 ~~  & ~~ 1 ~~  &  $\alpha_{opt}= 1.1704$  &  ~~ 1.4  ~~ & ~~   1.6 ~~   & ~~  1.8 ~~  \\
		\hline
		Iters   & 48            &  44         &  36   &  24  &   32  &  42  &  53   \\
		CPU  &  0.013      &  0.012   &   0.011      &  0.010       &   0.011      &   0.012      &   0.014    \\  
		\hline
	\end{tabular}
\end{table}   
 \end{example}

\begin{example}\label{Ex2} \rm
	In this example, we use four matrices TOLS340,TOLS1090, TOLS2000 and TOLS400 from Matrix Market  (\url{https://math.nist.gov/MatrixMarket/}). The Tolosa matrices are encountered in the stability analysis of an aircraft flight model. The notable modes of this system are characterized by complex eigenvalues, with their imaginary parts situated within a specified frequency range. The objective is to determine the eigenvalues that possess the largest imaginary components. This problem has been studied at CERFACS (Centre Européen de Recherche et de Formation Avancée en Calcul Scientifique) in collaboration with the Aerospatiale Aircraft division.	These matrices with their generic properties, including  the size ($n$), the number of nonzero entries ($nnz$) and the condition number estimation (Cond), are presented in Table \ref{Tbl2}. 
\begin{table}[!htp]
	\centering
	\caption{Matrix properties in Example \ref{Ex2}.  \label{Tbl2}}\vspace{0.25cm}
	\begin{tabular}{|c|ccc| } \hline
		matrix           &  ~~~~$n$~~~~    &  ~~~~ $nnz$ ~~~~  & ~~~~ Cond ~~~~ \\ \hline
		TOLS 340    &   340       &   2196        &        $2.4e+05$              \\ 
		TOLS1090   &   1090     &    3546       &        $2.1e+06$              \\ 
		TOLS2000   &   2000     &   5184        &        $6.9e+06$             \\ 
		TOLS4000   &   4000     &   8784        &        $2.7e+07$               \\ \hline		
	\end{tabular}
\end{table}   
We set  $A_1=$ TOLS 340, TOLS1090,   TOLS2000 and TOLS4000. In this case, we have $p=n=340,1090,2000$ and $4000$. For each system we set $q=n$ and $A_2=6I_n$, where $I_n$ is the identity matrix of order $n$. The entries of the vectors $b_1$ and $b_2$  are set to be all ones. 
 We solve the system \eqref{Eq06}  with restarted version of GMRES(10) in conjunction with the PBS preconditioner. Left-preconditioning is used.  
 A zero vector is used as an initial guess and the iteration is stopped as soon as the residual of the original system is reduced by a factor of $10^{11}$. The parameter $\alpha$ in the PBS method is set to be $1$.  {We also compare the numerical results of the PBS preconditioner with the preconditioners induced by  the block splittings \eqref{BS1}, \eqref{BS2}, \eqref{BS3}, and \eqref{BS4} (denoted by  BS1, BS2, BS3 and  BS3) applied to the system \eqref{Eq06}.  In the implementation of the preconditioners  $\mathcal{M}_i$ $(i=1,2,3,4)$ within  the GMRES method, we need to solve systems of the form  $\mathcal{M}_i (z_1; z_2;z_3)=  (r_1; r_2;r_3) $ which can be accomplished as follows (see \cite{LiAMC,Xin}). }
 
 \bigskip
    
  {\textbf{Algorithm 2. } Solution of  $\mathcal{M}_{1}(z_1;z_2;z_3)=(r_1;r_2;r_3)$ \\[-0.8cm]
 	\begin{itemize}
 		\item [1.] Set $z_1=r_1$;\\ [-0.8cm]
 		\item [2.] Solve $Pz_2=r_2$ for $z_2$;\\ [-0.8cm]
 		\item [3.] Set $z_3=r_3$.
 \end{itemize} }
 
  {\textbf{Algorithm 3. } Solution of  $\mathcal{M}_{2}(z_1;z_2;z_3)=(r_1;r_2;r_3)$ \\[-0.8cm]
	\begin{itemize}
		\item [1.] Set $z_1=r_1$;\\ [-0.8cm]
		\item [2.] Set $z_3=r_3$;\\ [-0.8cm]
		\item [3.] Solve $Pz_2=r_2-A_2^Tz_3$ for $z_2$.
\end{itemize} }
 
  {\textbf{Algorithm 4. } Solution of  $\mathcal{M}_{3}(z_1;z_2;z_3)=(r_1;r_2;r_3)$ \\[-0.8cm]
 	\begin{itemize}
 		\item [1.] Solve $Pz_2=r_2$ for $z_2$;\\ [-0.8cm]
 		\item [2.] Compute $z_1=r_1-A_1 z_2$;\\ [-0.8cm]
 		\item [3.] Set  $z_3=r_3$.
 \end{itemize} }

  {\textbf{Algorithm 5. } Solution of  $\mathcal{M}_{4}(z_1;z_2;z_3)=(r_1;r_2;r_3)$ \\[-0.8cm]
	\begin{itemize}
		\item [1.] Set $z_3=r_3$\\ [-0.8cm]
		\item [2.] Solve $Pz_2=r_2-A_2^T z_3$ for $z_2$;\\ [-0.8cm]
		\item [3.] Set  $z_1=r_1-A_1z_2$.
\end{itemize} }

  Numerical results have been presented in Table \ref{Tbl2}. In this table the number of iterations (Iter), elapsed CPU time in second (CPU) ,
 $$Rel=\frac{\|\textbf{r}_k\|_2}{\|\textbf{r}_0\|_2},$$
with $r_k$ is the residual vector at iteration $k$,  
 $$Err=\frac{\|\textbf{x}_k-\textbf{x}^*\|_2}{\|\textbf{x}^*\|_2},$$
where $\textbf{x}_{k}$ is the computed solution at iteration $k$  and $\textbf{x}^*$ is the exact solution to the system.
To show the effectiveness of the preconditioners we also examined the restarted version of GMRES(10) method without preconditioning (No-Prec). A dagger ($\dag$) shows that the method has not converged in 1000 iterations. As the numerical results show all the preconditioners are effective. We also see that the PBS preconditioner is more efficient that the other preconditioners in terms of the number of iterations and the CPU time.     
	  
\end{example}

\begin{table}[!htbp]
	\centering
	\caption{Numerical results for Example \ref{Ex2}. \label{Tbl3}} \vspace{0.5cm}
	\begin{tabular}{|c|l|c|c|c|c|c|c|}
		\hline
		~~~$m\times n$                &               & PBS          &   BS1         &  BS2            & BS3       & BS4       & No-Prec\\  \hline
		$640\times 340$        & Iter         &  14            &  32             &  18               & 34                 &  16 &$\dag$ \\ 
		                                     & CPU      & 0.12          &  0.23          &  0.19           & 0.23              &  0.18 &\\
		                                     & RES       & 9.13e-12  & 3.79e-12  & 3.52e-13    & 7.31e-12     &  3.21e-12  &\\
		                                     & ERR      & 9.22e-12  & 3.07e-11  & 9.16e-13    & 9.06e-12     &  6.39e-12  & \\   \hline
		$2180 \times 1090$   & Iter         &  14            &  32             &  16               & 34                &  15  &$\dag$ \\ 
                                       		 & CPU      & 0.15          &  0.21          &  0.18           & 0.23              & 0.16  & \\
		                                     & RES       & 2.58e-12  & 2.71e-12  & 6.05e-12    & 6.26e-12      & 6.43e-12  &\\
		                                     & ERR      & 3.61e-12  & 1.83e-11  & 1.76e-11    & 7.87e-12      & 2.57e-11 &\\   \hline

		$4000 \times 2000 $   & Iter         &  14            &  32             &  16               & 34               & 15 & $\dag$ \\ 
                                              & CPU      & 0.17          &  0.24          &  0.19           & 0.25             & 0.19 & \\
		                                      & RES       & 1.91e-12  & 4.12e-12  & 4.42e-12    & 4.78e-12    & 4.5e-12 &\\
		                                      & ERR      & 3.02e-12  & 3.45e-11  & 1.28e-11    & 5.05e-12    & 1.80e-11 &\\   \hline

		$8000 \times 4000$    & Iter         &  13            &  30             &  16               & 34               & 15  & $\dag$ \\ 
		                                      & CPU      & 0.21          &  0.34          &  0.25           & 0.40             & 0.24  & \\
		                                      & RES       & 7.96e-12  & 8.98e-12  & 3.04e-12    & 3.07e-12    & 3.6e-12  & \\
	                                          & ERR      & 1.81e-11  & 1.44e-10  & 9.70e-12    & 5.23e-12    & 1.29e-11  & \\   \hline

		\end{tabular}
\end{table}

\begin{example}\label{Ex3} \rm
We consider the partial differential equations
\begin{equation} 	  
	\left\{\begin{array}{rll}
		-\Delta u +\sin(x+y)\frac{\partial u}{\partial x} +\cos(x-y) \frac{\partial u}{\partial y}  + 50(x+y) u = f, & {\rm in} & \Omega=[0,1]\times [0,1], \\
		u = g, & {\rm on} & \partial\Omega,
	\end{array}
	\right.	
\end{equation}
where the function $g$ is known.The discretization of this equation employing the second-order central difference scheme on an $(n_0 + 2)$ grid in all directions of  $\Omega$, with a mesh size of  $h = 1/(n_0+1)$, yields a system of linear equations with the coefficient matrix $A_1  \in \mathbb{R}^{n \times n},$  in which  $n=n_0^2$. We set $A_2=0.7 I_n$, where $I_n$ is the identity matrix of order $n$.
  {We set $n_0=85,90,95,100,105$ and $110$. Consequently, the resulting matrices $A_1$ are of orders $7225$, $8100$,  $9025$, $10000$, $11025$ and $12100$, respectively.  In this case, the order of matrices $A$ is  $14450\times 7225$, $16200\times 8100$, $18050\times 9025$, $20000\times 10000$, $22050\times 11025$ and $24200\times 12100$.}
All the assumptions are as the previous example, except here we have used the complete version of GMRES. All the other assumptions are the same as in the previous example. Numerical results are presented in Table \ref{Tbl3}.  As we observe the PBS preconditioner outperforms the other preconditioners.   
 
\begin{table}[!htbp]
	\centering
	\caption{Numerical results for Example \ref{Ex3}. \label{Tabl4}} \vspace{0.5cm}
	\begin{tabular}{|c|c|c|c|c|c|c|c|c|}
		\hline
	$n_0$	&~~~$m\times n$       &               & PBS          &   BS1         &  BS2            & BS3  & BS4            & No-Prec\\  \hline
 85 & $ 14450\times 7225$       
& Iter  &  4       &  9          &  6             & 9         & 5         &   $\dag$\\ 
&& CPU   & 0.49     &  0.81       &  0.59          & 0.74      & 0.53      & \\
&& RES   & 7.49e-12 & 5.86e-12    & 4.63e-13       & 1.47e-12  & 1.43e-12  &\\
&& ERR   & 4.30e-09 & 1.01e-09    & 1.68e-09       & 1.68e-09  & 1.64e-9   &\\   \hline
90 &$16200 \times 8100$    
& Iter &  4        &  9          &  6             & 9         & 5         & $\dag$  \\ 
&& CPU  & 0.57      &  0.96       &  0.69          & 0.86      & 0.62      &         \\
&& RES  & 7.52e-12  & 5.87e-12    & 5.22e-13       & 1.65e-12  & 1.56e-12  &         \\
&& ERR  & 3.43e-09  & 1.88e-09    & 1.88e-09       & 1.88e-09  & 1.28e-9   &         \\   \hline
95 &$18050\times 9025 $    
&Iter  &  4        &  9        &  6        & 9            & 5         & $\dag$ \\ 
&&CPU   & 0.67      &  1.10     &  0.80     & 1.01         & 0.72      & \\
&&RES   & 7.53e-12  & 5.88e-12  & 5.71e-13  & 1.75e-12     & 1.75e-12  & \\
&&ERR   & 5.85e-09  & 1.65e-09  & 1.64e-09  & 1.65e-9      & 2.16e-9   &\\   \hline

100 &$20000\times 10000 $    
&Iter  &  4        &  9        &  6        & 9            & 5         & $\dag$ \\ 
&&CPU   & 0.78      &  1.28     &  0.94     & 1.18         & 0.86      & \\
&&RES   & 7.55e-12  & 5.91e-12  & 6.67e-13  & 1.90e-12     & 1.94e-12  & \\
&&ERR   & 4.87e-09  & 2.29e-09  & 2.29e-09  & 2.29e-9      & 2.29e-9   &\\   \hline

105 &$22050\times 11025 $    
&Iter  &  4        &  9        &  6        & 9            & 5         & $\dag$ \\ 
&&CPU   & 0.90      &  1.46     &  1.09     & 1.39         & 0.98      & \\
&&RES   & 7.56e-12  & 5.91e-12  & 7.28e-13  & 2.09e-12     & 2.14e-12  & \\
&&ERR   & 4.69e-09  & 2.08e-09  & 2.08e-09  & 2.08e-9      & 2.08e-9   &\\   \hline

110 &$24200\times 12100 $    
&Iter  &  4        &  9        &  6        & 9            & 5         & $\dag$ \\ 
&&CPU   & 1.03      &  1.66     &  1.25     & 1.58         & 1.12      & \\
&&RES   & 7.57e-12  & 5.93e-12  & 8.09e-13  & 2.28e-12     & 2.31e-12  & \\
&&ERR   & 7.69e-09  & 3.19e-09  & 3.19e-09  & 3.19e-9      & 3.19e-9   &\\   \hline
		
	\end{tabular}
\end{table}
	
\end{example}
\section{Conclusions}\label{Sec5}
We have presented a parameterized block splitting iteration method for solving  indefinite least squares problem.  {We have proved that the iteration method is conditionally convergent.} Optimal value of the parameter has been computed. Eigenpair properties of the preconditioned matrix has been studied. We have compared the numerical results of the preconditioner induced by the proposed method with three recently presented preconditioners. The numerical results show that the proposed preconditioner in more efficient that others.  

\section*{Acknowledgements}
 {The author sincerely thanks the reviewers for their meticulous review of the paper and for their constructive comments and suggestions}   

\section*{Conflicts of interest}
This work does not have any conflicts of interest.


\end{document}